\documentclass[a4paper,12pt]{amsart}

\usepackage{enumerate, booktabs, mathtools, mathrsfs}
\usepackage{amsmath, amsfonts, amssymb, amsthm}
\usepackage{mathtools} 
\usepackage{color}
\usepackage[dvipsnames]{xcolor}
\usepackage{relsize}
\usepackage[all]{xy}
\usepackage[height=22.6cm, width=15.5cm]{geometry}

\usepackage{needspace} 

\providecommand{\bysame}{\leqavevmode\hbox to3em{\hrulefill}\thinspace}

\numberwithin{equation}{section}

\theoremstyle{plain}
\newtheorem{nnthm}{Theorem}
\newtheorem{thm}{Theorem}[section]
\newtheorem{prop}[thm]{Proposition}
\newtheorem{lem}[thm]{Lemma}

\newtheorem*{Prop3.2}{Proposition 3.2}
\newtheorem*{thm4.2}{Theorem 4.2}

\newtheorem*{thm*}{Theorem}

\theoremstyle{definition}

\theoremstyle{remark}
\newtheorem*{rem}{Remark}
\newtheorem*{ex}{Example}

\newcommand{\mbb}[1]{\mathbb{#1}}
\newcommand{\wt}[1]{\widetilde{#1}}

\newcommand{\ol}[1]{\overline{#1}}

\newcommand{\abs}[1]{\lvert #1\rvert}

\renewcommand{\leq}{\leqslant}
\renewcommand{\geq}{\geqslant}

\def\slsf{\slshape \sffamily }

\DeclareMathOperator{\id}{id}

\DeclareMathOperator{\Aut}{Aut}

\DeclareMathOperator{\rk}{rk}

\DeclareMathOperator{\eps}{\varepsilon}


\title[The Levi\, problem\, over\, generalized\, Hirzebruch manifolds]%
{The Levi\, problem\, over\, generalized \\[3pt] Hirzebruch manifolds}

\author{S. Ivashkovych}
\address{D\'epartement de Math\'ematiques, Universit\'e de Lille, Villeneuve
d'Ascq, France}
\email{serge.ivashkovych@univ-lille.fr}

\author{C. Miebach}
\address{Univ.~Littoral C\^ote d'Opale, UR 2597, LMPA, Laboratoire de 
Math\'ematiques Pures et Appliqu\'ees Joseph Liouville, F-62100 Calais, France}
\email{christian.miebach@univ-littoral.fr}

\author{V. Shevchishin}
\address{Faculty of Mathematics and Computer Sciences, 
University of Warmia and Mazury, Olsztyn, Poland}
\email{vsevolod@matman.uwm.edu.pl}

\dedicatory{In memoriam Alan T.~Huckleberry}

\begin{document}

\begin{abstract}
We review classical methods to solve the Levi problem in the presence of
symmetries, established by Hirschowitz and by Grauert-Remmert-Ueda. We then
illustrate these methods by solving the Levi problem in some new situations,
namely generalized Hirzebruch manifolds and 
primary Hopf surfaces of non-diagonal type.
\end{abstract}

\maketitle

\section{Introduction}

One of the central problems in the theory of several complex variables is the
so-called \emph{Levi problem}, which can be stated as follows. Let $X$ be a
complex manifold and let $(D,p)$ be a locally Stein domain over $X$. Under which
conditions is $D$ Stein?
For $X=\mbb{C}^n$ the Levi problem was first solved by Oka~\cite{O} who showed
that every locally Stein domain over $\mbb{C}^n$ is Stein. This fundamental
result was extended by Docquier and Grauert to locally Stein domains over
arbitrary Stein manifolds, see~\cite{DG}. In~\cite{Fj} Fujita proved that the
only locally Stein domain over $\mbb{P}^n$ that is not Stein is
$(\mbb{P}^n,\id)$.

\smallskip In this note we are going to present two methods that allow to
 solve the Levi
problem in the presence of symmetries. The first one was developed by Ueda
in~\cite{U} in order to extend Fujita's theorem to locally Stein domains over
Grassmann manifolds. His main insight was that one can make use of
Matsushima's and Morimoto's work~\cite{MM} on holomorphic principal bundles of
complex reductive groups together with deep results by Grauert and
Remmert from~\cite{GR1}. This is related to~\cite{Nem}.

\smallskip We will illustrate these ideas by solving the Levi problem for 
locally Stein
domains over generalized Hirzebruch manifolds $X_k$, i.e., for the total
spaces of the projectivization of the holomorphic vector bundle $\mathcal{O}
\oplus \mathcal{O}(k)\to {\rm{Gr}}_d(\mbb{C}^n)$ for some $k\geq0$ where
${\rm{Gr}}_d(\mbb{C}^n)$ is the Grassmann manifold of $d$-dimensional linear
subspaces of $\mbb{C}^n$. In particular, this includes the solution of the Levi
problem for domains over Hirzebruch surfaces.

\begin{nnthm}
\label{Thm:Hirzebruch}
Let $D$ be a locally Stein domain over a generalized Hirzebruch manifold
$X_k$ for $k\geq 0$. 
If $D$ is not Stein, then one of the following possibilities occurs:
\begin{enumerate}[(1)]
\item $D=q_k^{-1}(V)$ for some locally Stein domain $V$ over ${\rm{Gr}}_d(\mbb{C}^n)$ (this
includes the case $D=X_k$);
\item $D$ is a finite unramified cover over a $1$-complete neighbourhood of the
exceptional divisor $E_\infty$;
\item $D$ is an finite unramified cover over $B\setminus E_\infty$, where $B$ is a
$1$-complete neighbourhood of the exceptional divisor $E_\infty$;
\item $D$ is a finite cover over $X_k\setminus E_\infty$ (this includes the case
$D=X_k\setminus E_\infty$);
\end{enumerate}
\end{nnthm}

Here, by saying that $D$ (or $B$) is $1$-complete we mean that after
contracting the exceptional divisor $E_\infty$, it becomes Stein.

As a second application, we will study locally Stein domains in certain
products of compact Riemann surfaces.

\begin{nnthm}
\label{Thm:Products}
Let $\Sigma_g$ be a compact Riemann surface of genus $g\geq0$ and let $D$ be a
locally Stein domain in $\Sigma_g\times \mbb{P}^1$. If $D$ is not Stein, then 
only the following cases are possible:
\begin{enumerate}[(1)]
\item $D= D_1\times \mbb{P}^1$, where $D_1$ is a domain in $\Sigma_g$ or
\item $D = \Sigma\times D_2$, where $D_2$ is a domain in $\mbb{P}^1$.
\end{enumerate}
\end{nnthm}

This result uses~\cite{GR1} together with a theorem of Brun, see~\cite{Brn},
that solves the Levi problem for domains in holomorphic fibre bundles over
Stein manifolds with fibre a compact Riemann surface.

\smallskip The second one of the aforementioned geometric methods was developed by
Hirscho\-witz in~\cite{Hir1} and \cite{Hir2} where he studied the Levi
problem for domains $(D,p)$ over arbitrary infinitesimally homogeneous manifolds
$X$ under the more restrictive assumption that $D$ admits a continuous
plurisubharmonic exhaustion function. Such domains will be called 
\emph{pseudoconvex}
in order to distinguish them from locally Stein domains.
We will illustrate Hirschowitz' techniques by proving the following result.

\begin{nnthm}
\label{Thm:Hopf}
Let $X$ be a primary Hopf surface of non-diagonal type. Then every pseudoconvex
domain $D\subsetneq X$ is Stein.
\end{nnthm}

This result completes the work begun in~\cite{LY} and~\cite{Mie}, 
where the case of diagonal primary Hopf surfaces was considered.
 
\section{The Levi problem over generalized Hirzebruch manifolds}

In this section we closely follow Ueda's arguments from~\cite{U} in order to
solve the Levi problem for domains over projectivized line bundles over
Grassmannians.

\subsection{Ueda's generalization of methods by Grauert and Remmert}

Recall that a Riemann domain over a complex manifold  $X$ is a pair $(D, p )$
where $D$ is a connected Hausdorff topological space and $p \colon D\to X$ is
local homeomorphism. This homeomorphism induces an obvious complex structure on
$D$ and $p$ becomes a local biholomorphism, see \cite{GR2}. The Riemann domain
$(D,p)$  is called \emph{locally Stein over the point} $x\in X$ if there exists
a Stein neighbourhood $U\ni x$ such that all connected components of $p^{-1}(U)$
are Stein. Moreover, $(D,p)$ is called locally Stein over $X$ if it is locally
Stein over every point of $X$.

\smallskip
For later use we note an important result by Grauert and Remmert which was
later generalized by Ueda, see~\cite[Satz~6]{GR1} and \cite[Section~4]{U}.
For its statement, recall that for a Riemann domain $(D,p)$ over a complex
manifold $X$ one can define the boundary $\partial D$,
see~\cite[Definition~4]{GR1}. Moreover, $\ol{D}:=D\cup \partial D$ can be given
the structure of a Hausdorff topological space such that $p\colon D\to X$
extends to a continuous map $\ol{p}\colon\ol{D}\to X$. Let $A\subset X$ be a
nowhere dense analytic set. A boundary point $x\in\partial D$ is called
\emph{removable along $A$} if it possesses a neighbourhood $U$ in $\ol{D}$ such
that $(U,\ol{p}|_U)$ is a schlicht domain satisfying $U\cap\partial D\subset
\ol{p}^{-1}(A)$. We denote by $\partial^rD$ the set of all boundary points
which are removable along $A$. It can be shown that
\begin{equation*}
(D^*:=D\cup\partial^rD,p^*:=\ol{p}|_{D^*})
\end{equation*}
is again a Riemann domain over $X$, called the \emph{extension of $D$ along
$A$}.

\begin{thm}[Grauert-Remmert, Ueda]\label{Thm:GRU}
Let $X$ be a complex manifold and let $(D,p)$ be a Riemann domain over $X$. Let
$A\subset X$ be an analytic subset of positive codimension such that $D$ is
locally Stein over $X\setminus A$.
\begin{enumerate}[(1)]
\item If no boundary point of $D$ is removable along $A$, then $D$ is locally
Stein over $X$.
\item The extension $p^*\colon D^*\to X$ along $A$ is locally Stein over $X$.
\end{enumerate}
\end{thm}

\subsection{Generalized Hirzebruch manifolds}

Let $M={\rm{Gr}}_d(\mbb{C}^n)$ be the Grassmann manifold of all $d$-dimensional
complex linear subspaces of $\mbb{C}^n$. Since the Picard group of $M$ is
isomorphic to $\mbb{Z}$, for every integer $k\geq0$ we have a holomorphic line
bundle $\mathcal{O}(k)\to M$, which we can compactify to a holomorphic
$\mbb{P}^1$-bundle $q_k\colon X_k\to M$. In other words, $X_k$ is the
projectivization of the holomorhic vector bundle
$\mathcal{O}\oplus\mathcal{O}(k)$ of rank $2$ over ${\rm{Gr}}_d(\mbb{C}^n)$.
Note that, since the line bundle $\mathcal{O}(k)$ is positive, the infinity
section $E_\infty$ in $X_k$ is an exceptional divisor. Hence, after contracting
$E_\infty$ the open set $X_k\setminus E_0$ (where $E_0$ is the zero section in
$X_k$) becomes a Stein space with an isolated singularity.

\begin{ex}
If $M=\mbb{P}^1={\rm{Gr}}_1(\mbb{C}^2)$, then $X_k$ is the $k$-th Hirzebruch
surface.
\end{ex}

In order to apply Ueda's ideas we will describe $X_k$ as the base of a
holomorphic principal bundle. For this, recall first that every $g\in
{\rm{GL}}(d,\mbb{C})$ can be uniquely written as $g=\lambda_g\wt{g}$ where
$\lambda_g\in\mbb{C}^*$ and $\wt{g}\in{\rm{SL}}(d,\mbb{C})$. The group
$G:={\rm{GL}}(d,\mbb{C})\times\mbb{C}^*$ acts on $V:=\mbb{C}^{n\times
d}\times\mbb{C}^2$ by
\begin{equation*}
(g,\mu)\cdot(Z,w):=\bigl(Zg^{-1},(\mu w_1,\lambda_g^k\mu w_2)\bigr).
\end{equation*}

\begin{rem}
Note that $G$ is complex reductive as it is the complexification of its maximal
compact subgroup ${\rm{U}}(d)\times{\rm{U}}(1)$.
\end{rem}

Let $\mbb{C}^{n\times d}_{d-1}$ be the analytic set of $Z\in\mbb{C}^{n\times d}$ 
such
that $\rk(Z)\leq d-1$ and let us consider the two analytic subsets
$A_1:=\mbb{C}^{n\times d}\times\{0\}$ and $A_2:=\mbb{C}^{n\times d}_{d-1}
\times\mbb{C}^2$ of $V$. One checks directly that $G$ acts freely and properly
on $\Omega:=V\setminus(A_1\cup A_2)$ as well as that $\Omega/G\cong X_k$. In
other words, we have a holomorphic $G$-principal bundle $\pi\colon\Omega\to
X_k$.

\begin{rem}
The exceptional divisor $E_\infty$ of $X_k$ coincides with $\pi\bigl(\{w_1=0\}
\bigr)$.
\end{rem}

For later use we note the following

\begin{lem}\label{Lem:OrbitClosure}
For $(Z,w)\in\Omega$ let $\ol{G\cdot(Z,w)}$ be the topological closure in $V$ of
the corresponding $G$-orbit. Then, we have
\begin{equation*}
\ol{G\cdot(Z,w)}=\ol{{\rm{GL}}(d,\mbb{C})\cdot Z}\times\{0\}\cup\{0\}\times
(\mbb{C}e_1)\cup \{0\}\times(\mbb{C}e_2).
\end{equation*}
In particular, the origin of $V$ lies in the closure of every $G$-orbit in
$\Omega$.
\end{lem}

\subsection{Locally Stein domains over $X_k$}

For every Riemann domain $(D,p)$ over $X_k$ we can define its pull-back to
$\Omega$ by
\begin{equation*}
\wt{D}:=\bigl\{(v,x)\in\Omega\times D;\ \pi(v)=p(x)\bigr\}
\end{equation*}
and hence obtain the commutative diagram
\begin{equation*}
\xymatrix{
\wt{D}\ar[r]^{\wt{\pi}}\ar[d]_{\wt{p}} & D\ar[d]^p\\
V\supset \Omega\ar[r]_\pi & X_k.
}
\end{equation*}
Here, $\wt{\pi}\colon\wt{D}\to D$ is again a $G$-principal bundle and we may
view $(\wt{D},\wt{p})$ as a Riemann domain over $V$.

We will need the following result of Matsushima and Morimoto,
see~\cite[Th\'eor\`emes 4 et 5]{MM}.

\begin{thm}\label{Thm:MM}
Let $G$ be a complex reductive group and let $\pi\colon P\to B$ be a holomorhic
$G$-principal bundle. Then, the total space $P$ is Stein if and only if the base
$B$ is Stein.
\end{thm}

We are now in position to prove Theorem~\ref{Thm:Hirzebruch}. For this, let
$(D,p)$ be a locally Stein domain over the generalized Hirzebruch manifold
$X_k$. It follows from Theorem~\ref{Thm:MM} that $\wt{p}\colon\wt{D}\to V$ is
locally Stein over $\Omega$.

Let $\partial^r\wt{D}$ be the set of boundary points of $\wt{D}$ which are
removable along $A_1\cup A_2$. If $\partial^r\wt{D}=\emptyset$, then it follows
from Theorem~\ref{Thm:GRU} that $\wt{D}$ is locally Stein over $V$. Then,
$\wt{D}$ is Stein and hence Theorem~\ref{Thm:MM} implies that $D$ is Stein as
well. We will therefore assume in the following that $\partial^r\wt{D}$ is not
empty. Let $(\wt{D}^*,\wt{p}^*)$ be the extension of $(\wt{D},\wt{p})$ along
$A_1\cup A_2$.

If $x\in\partial^r\wt{D}$, then $\wt{p}^*(x)\in\ol{G\cdot(Z,w)}\subset A_1\cup
A_2$ for some $(Z,w)\in\Omega$. Using Lemma~\ref{Lem:OrbitClosure} we will
distinguish several cases.

\paragraph{\slsf Case 1}

If $\partial^r\wt{D}$ contains a point lying over the origin of $V$, then
$\wt{D}^*$ contains a schlicht domain biholomorphic to a neighborhood of $0\in
V$. Since every $G$-orbit in $\Omega$ contains $0$ in its closure, we obtain
$D=X_k$ in this case.

\paragraph{\slsf Case 2}

Suppose that $\partial^r\wt{D}$ contains a point lying over $(Z,0)$ for some
$0\not=Z\in\mbb{C}^{n\times d}_{d-1}$. Using again Lemma~\ref{Lem:OrbitClosure}
we see that in this case $\wt{D}^*$ contains a schlicht domain biholomorphic to
$\{Z\}\times\mbb{C}^2$. Hence, $D$ contains with every point the whole
$\mbb{P}^1$-fiber passing through this point. This means that $D$ has the form
$q_k^{-1}(V)$ for some domain $V$ over ${\rm{Gr}}_d(\mbb{C}^n)$.

\paragraph{\slsf Case 3}

Suppose that $\partial^r\wt{D}$ contains a point lying over $(0,w)$ with
$w\not=0$. Then $\wt{D}^*$ contains a schlicht neighbourhood over $U\times
\Delta^2_{\eps}(w^*)$ for some $\eps>0$ where $U=U_0\setminus\mbb{C}^{n\times
d}_{d-1}$ for some neigborhood $U_0$ of $0\in\mbb{C}^{n\times d}$. We will
always assume that $U_0$ is invariant under multiplication by $\lambda_g$ with
$\abs{\lambda_g}<1$.

Acting on this set by $g\in {\rm{GL}}(d,\mbb{C})$ with $\abs{\lambda_g}<1$ we
get that $\wt{D}^*$ contains a domain spread over
\begin{equation}\label{cone1}
U\times \Delta_{\eps}(w_1^*) \times \big\{\,|w_2| \geq |w_2^*|\,\big\}.
\end{equation}
Indeed, take an arbitrary $Z\in U$ and notice that $\lambda_gZ\in U$
for any $\abs{\lambda_g}<1$. Then act by $\lambda_g^{-1}$ to get a point in
$\wt{D}^*$ over the point $x=(Z,w_1,\lambda^{-k}w_2)\in\Omega$. This means
that $\wt{D}^*$ may contain $d\leq k$ points over this $x$ and they satisfy
\eqref{cone1}.

\paragraph{\slsf Subcase 3a}

Let  $w_1^* \not =0$. Acting now by $\mu$ we obtain that $\wt{D}^*$ contains a
domain spread over the  cone
\begin{equation}\label{cone2}
U\times 
\bigg(\bigg\{|w_2| \geq \frac{|w_2^*|}{|w_1^*|}\,|w_1|  \bigg\}
\setminus \big\{w_1=0\big\}\bigg).
\end{equation}
Indeed, for an arbitrary $w_1\not=0$ find $\mu$ such that  $w_1 = \mu w_1^*$
and act on $(w_1^*, w_2)$, which satisfies \eqref{cone1}, by this $\mu$ to
get $(w_1,\mu w_2)$. Now, since we had $|w_2|> |w_2^*|$ we will obtain that
$|\mu w_2|> |\mu w_2| = |w_2^*|\frac{|w_1|}{|w_1^*|}$ as required. Therefore, in
this case $D$ is finitely spread over $B\setminus E$ for some $1$-complete
domain $B$ containing  $E$. I.e., the case (3)  occurs if $w_2^*\not=0$ and
the case (4) if $w_2^*=0$. Since $\wt{D}^* \to U$ is an unramified and the
fundamental group of $H_k\setminus E$ is $\mbb{Z}_k$ for $k\geq 2$ we see that
covering number $d$ should be a divisor of $k$.

\paragraph{\slsf Subcase 3b}

Now let $w_1^*=0\not= w_2^*$. Then $\wt{D}^*$ contains a schlicht neighbourhood
over $U\times \Delta_{\eps}(0)\times \Delta_{\eps}(w_2^*)$ for some $\eps >0$.
Acting by $\mu$ we obtain that $\wt{D}^*$ contains a domain spread over the cone
\begin{equation}\label{cone3}
U\times \bigg\{|w_2| \geq \frac{|w_2^*|}{\eps }\;|w_1|\bigg\}.
\end{equation}
Indeed, for $|w_1|>\eps$ take $\mu = \frac{\eps}{w_1}$ and acting with this
$\mu$ get from $(w_1,w_2)$ the point $(\eps , \frac{\eps w_2}{w_1})$. Then due
to \eqref{cone1} we have
\begin{equation*}
\left|\frac{\eps w_2}{w_1}\right|\geq |w_2^*|,
\end{equation*}
which is \eqref{cone3}. In this case $D$ is a $1$-complete neighbourhood of the
exceptional divisor $E_\infty$, i.e., the case (2) occurs.

In the remaining part of this section we will prove Theorem~\ref{Thm:Products}. 

Let $D$ be a locally Stein domain in $\Sigma_g\times\mbb{P}^1$, where
$g\geq 0$. We denote by $\Omega$ the domain $\mbb{C}^2\setminus\{0\}$ and
consider the $\mbb{C}^*$-principal bundle $\pi\colon\Sigma_g\times\Omega\to
\Sigma_g\times\mbb{P}^1$. Let $\wt{D}$ be the inverse image of $D$ under $\pi$
and consider $\wt{D}$ as a domain in $\Sigma_g\times\mbb{C}^2$. According to
Theorem~\ref{Thm:MM} $\wt{D}$ is locally Stein at all points apart possibly of
the points in $A:=\Sigma_g\times \{0\}$. It is also $\mbb{C}^*$-invariant under
the action
\begin{equation*}
\lambda\cdot (s, w_1,w_2) = (s,\lambda w_1, \lambda w_2).
\end{equation*}
Denote by $\wt{D}^*$ the extension of $\wt{D}$ along $A$. Then, $\wt{D}^*$ is
locally Stein in $\Sigma_g\times\mbb{C}^2$ by Theorem~\ref{Thm:GRU}.

According to~\cite{Brn}, either $\wt{D}^*$ is Stein or $\wt{D}^*=\Sigma_g\times
U$ for a Stein domain $U\subset\mbb{C}^2$.

In the first case, either $D$ is Stein (if $\wt{D}^*$ does not contain a point
in $A$) or $D$ is of the form $D_1\times\mbb{P}^1$ for some domain $D_1\subset
\Sigma_g$. In the second case, $D$ is of the form $\Sigma_g\times D_2$ for some
domain $D_2\subset \mbb{P}^1$.

{\needspace{5\baselineskip} 
\section{The Levi problem in non-diagonal primary Hopf surfaces}%
In this section we are going to explain how Hirschowitz' methods lead to a
solution of the Levi problem for locally Stein domains in non-diagonal primary
Hopf surfaces. This completes the results of~\cite{LY} and~\cite{Mie} where
diagonal primary Hopf surfaces were considered.}

In the first subsection we are going to review Hirschowitz' techniques
from~\cite{Hir2} in a slightly more general form as presented in~\cite{Mie}.

\subsection{Hirschowitz' techniques}

Let $X$ be a complex manifold with holomorphic tangent bundle $TX\to X$, and let
$\pi\colon\mbb{P}TX\to X$ be the projectivized holomorphic tangent bundle. A
continuous function on $X$ is said to be \emph{strictly} pluri\-subharmonic on
$X$ if it is everywhere locally the sum of a continuous plurisubharmonic and a
smooth strictly plurisubharmonic function.

For every plurisubharmonic function $\varphi\in\mathscr{C}(X)$ we define
\begin{equation*}
C(\varphi):=\mbb{P}TX\setminus
\bigl\{[v]\in\mbb{P}TX;\text{ $\varphi$ is smooth around
$\pi[v]$ and } \partial\varphi(v)\not=0\bigr\}
\end{equation*}
as well as
\begin{equation}\label{Eqn:C(X)}
C(X):=\bigcap_{\substack{\varphi\in\mathscr{C}(X)\\\text{plurisubharmonic}}}
C(\varphi).
\end{equation}
Every set $C(\varphi)$ (and thus $C(X)$) is closed in $\mbb{P}TX$. The next
result is a slight gene\-ralization of~\cite[Proposition~1.5]{Hir2},
cf.~\cite[Lemma~2.3]{Mie}.

\begin{prop}\label{Prop:StrPsh}
Let $X$ be a complex manifold and let $\Omega:=X\setminus\pi\bigl(C(X)\bigr)$.
Then there is a plurisubharmonic function $\psi\in\mathscr{C}(X)$ which is
strictly plurisubharmonic on $\Omega$.
\end{prop}

In what follows we say that a complex manifold $X$ is \emph{pseudoconvex} if
there is a continuous plurisubharmonic exhaustion function $\rho\colon
X\to\mbb{R}^{>0}$.

The next result is \cite[Propositions~2.6 and 3.4]{Hir2}, see
also~\cite[Lemma~2.4]{Mie}.

\begin{prop}\label{Prop:InnerIntCurve}
Let $X$ be a pseudoconvex complex manifold and let $\gamma\colon U\to X$ be
an integral curve of a holomorphic vector field on $X$, where $U$ is a domain in
$\mbb{C}$. If $\gamma'(U)$ meets $C(X)$, then $\gamma'(U)$ is contained in
$C(X)$. If $X$ admits a smooth plurisubharmonic exhaustion function, then
$\gamma'(U)\subset C(X)$ implies that $\gamma(U)$ is relatively compact in $X$.
In particular, in this case we have $U=\mbb{C}$.
\end{prop}

Let us explain how these results imply Theorem~\ref{Thm:Products} for $g=0,1$.

Let $\Sigma_g$ be a compact Riemann surface of genus $g=0$ or $g=1$ and let
$X_g:=\Sigma_g\times\mbb{P}^1$. Since the surface $X_g$ is homogeneous under the
identity component $G=\Aut^0(X)$ of its automorphism group,
Hirschowitz' methods can be applied directly in order to describe locally
Stein domains $D$ in $X_g$.

First, let us note that every locally pseudoconvex domain $D\subset X_g$ admits
a continuous plurisubharmonic exhaustion function due to~\cite[Th\'eor\`eme~2.1,
Lemma~2.4]{Hir1}. Let us consider the set $C(D)\subset\mbb{P}TX$ defined as in
\eqref{Eqn:C(X)}. If $C(D)=\emptyset$, the domain $D$ is Stein due to
Proposition~\ref{Prop:StrPsh}. Let us therefore assume that $D$ is not Stein,
hence that $C(D)\cap \mbb{P}T_xX$ is non-empty for some $x\in D$.

The proof of Proposition~\ref{Prop:InnerIntCurve} shows that the isotropy group
$G_x$ acts on $C(D)\cap\mbb{P}T_xX$, see also~\cite[Proposition~1.6]{Hir2}. If
$g=0$, then $G={\rm{SL}}(2,\mbb{C})\times {\rm{SL}}(2,\mbb{C})$ and $G_x=B\times
B$ where $B$ is conjugate to the group of upper triangular matrices in
${\rm{SL}}(2,\mbb{C})$. If $g=1$, then $G=
\Sigma_1\times{\rm{SL}}(2,\mbb{C})$ and $G_x=\{0\}\times B$. In both cases $G_x$
acts on $\mbb{P}T_xX\cong\mbb{P}^1$ with two fixed points and one open orbit,
the two fixed points corresponding to the tangent lines of the horizontal and
vertical fibrations of $\Sigma_g\times\mbb{P}^1$. If $C(D)$ meets the open
$G_x$-orbit, then $D=X_g$ according to~\cite[Proposition~3.4]{Hir2} because all
the elements of $C(D)$ generate inner integral curves in $D$. If $C(D)$
coincides with one of the two $G_x$-fixed points, then $D$ contains with $x\in
D$ either the horizontal or the vertical fiber of $X_g=\Sigma_g\times\mbb{P}^1$.
Since neither of these fibers has a strongly pseudoconvex neighborhood, $C(D)$
must contain elements lying over points in $D\setminus\{x\}$ as well.
Consequently, we may iterate the above argument and conclude that
$D=D_1\times\mbb{P}^1$ with $D_1\subset\Sigma_g$ Stein or $D=\Sigma_g\times D_2$
with $D_2\subset\mbb{P}^1$ Stein or $D=X_g$.

\subsection{Locally Stein domain in non-diagonal Hopf surfaces}

Let us first recall the construction of a Hopf surface of non-diagonal type.
The map $f\colon\mbb{C}^2\setminus\{0\}\to\mbb{C}^2\setminus\{0\}$ defined by
the formula
\begin{equation*}
f(z,w)=(\lambda^kz+w^k,\lambda w)\quad\bigl(k\geq1,0<\abs{\lambda}<1\bigr)
\end{equation*}
generates a free proper $\mbb{Z}$-action on $\mbb{C}^2\setminus\{0\}$. The
corresponding quotient $X:=(\mbb{C}^2\setminus\{0\})/\mbb{Z}$ is a primary Hopf
surface of non-diagonal type. We denote the quotient map
$p\colon\mbb{C}^2\setminus\{0\}\to X$ by $p(z,w)=[z,w]$. The group
$G:=\mbb{C}^*\times\mbb{C}$ acts on $X$ by
\begin{equation*}
(a,b)\cdot[z,w]:=\bigl[a^k(z+bw^k),aw\bigr].
\end{equation*}
One can show that $\Aut(X)\cong G/\langle f\rangle$.
The group $G$ acts on $X$ with two orbits. The orbit $E:=G\cdot[1,0]$ is an
elliptic curve in $X$ and the orbit $X^*:=G\cdot[0,1]\cong G/\langle
f\rangle\cong \mbb{C}^*\times\mbb{C}^*$ is open in $X$.

\smallskip 
Let $D\subsetneq X$ be a locally Stein domain. We shall show that $D$ is
Stein. We first note that $D$ cannot contain $E$. This follows from the fact
that, since $X^*$ is Stein, the domain $D^*:=D\cap X^*$ is Stein as well, hence
must have connected boundary. Using the $G$-action on $X$ and the fact that $D$
is not $G$-invariant, we obtain a continuous plurisubharmonic exhaustion
function of $D$ like in the proof of~\cite[Proposition~4.1]{Mie}.

In the next step we shall show that this plurisubharmonic exhaustion
function of $D$ is strictly plurisubharmonic on $D^*$. According
to Propositions~\ref{Prop:StrPsh} and~\ref{Prop:InnerIntCurve} we must only
exclude the possibility that there exists a one-parameter group of $A\subset G$
such that $A\cdot x$ is relatively compact in $D$ for some $x\in D^*$. Let
$A\subset G=\mbb{C}^*\times\mbb{C}$ be such a one-parameter group. If the
projection from $A$ to $\mbb{C}^*$ is non-trivial, then $A$ acts transitively on
$E$. Therefore the closure of $A\cdot x$ must contain $E$. Otherwise we have
$A\cong\mbb{C}$ acting by $t\cdot[z,w]=[z+tw^k,w]$. Without loss of generality
we may suppose that $x=[0,1]\in D$. Set $t_m:=t_0\lambda^{-m}$ and note that
$\lim_{m\to\infty}t_m=\infty$. Then
\begin{equation*}
t_m\cdot[0,1]=[t_m,1]=[\lambda^mt_m+m\lambda^{(m-1)k},\lambda^m]=
[t_0+m\lambda^{(m-1)k},\lambda^m]
\end{equation*}
for all $m\in\mbb{Z}$. For $m\to\infty$ this sequence converges to $[t_0,0]$.
Thus $E$ is contained in the closure of $A\cdot[0,1]$. We conclude that
existence of $A$ implies that $D$ contains $E$, which is not possible.

Since $D$ has a plurisubharmonic exhaustion that is strictly plurisubharmonic
on $D^*$ and since $D\cap E$ is Stein, an application of~\cite[Lemma~5.3]{Mie}
proves that $D$ is Stein.

\end{document}